\documentclass[11pt,leqno,a4paper]{article}

\usepackage{amsmath, amssymb, amsthm, amsfonts}
\usepackage{mathrsfs}
\usepackage{color}
\usepackage{hyperref}
\theoremstyle{plain}
\newtheorem{thm}{Theorem}

\newtheorem{lem}{Lemma}

\theoremstyle{remark}
\newtheorem{rem}{Remark}

\renewcommand{\Re}{{\rm Re\,}}

\def\K{\mathop{\mbox{\bf\Large K}}}

\numberwithin{equation}{section}

\begin{document}

\title{Continued fraction expression of the Mathieu series}
\author{Xiaodong Cao*, Yoshio Tanigawa and Wenguang Zhai }
\date{}

\maketitle
\footnote[0]{* Corresponding author.}

\footnote[0]{\emph{E-mail address}: caoxiaodong@bipt.edu.cn~(X.D. Cao), tanigawa@math.nagoya-u.ac.jp~(Y. Yanigawa), zhaiwg@hotmail.com~(W.G. Zhai)}



\footnote[0]{This work is supported by
the National Natural Science Foundation of China (Grant No.11171344) and the Natural
Science Foundation of Beijing (Grant No.1112010).}

\footnote[0]{
Xiaodong Cao:
Department of Mathematics and Physics,
Beijing Institute of Petro-Chemical Technology,
Beijing, 102617, P. R. China \\
Yoshio Tanigawa:
Graduate School of Mathematics,
Nagoya University,
Nagoya, 464-8602, Japan\\
Wenguang Zhai:
Department of Mathematics,
China University of Mining and Technology,
Beijing 100083, P. R. China}

\begin{abstract}
In this paper, we represent a continued fraction expression of Mathieu series by a continued fraction formula of Ramanujan. As application, we obtain some new bounds for Mathieu series.
\end{abstract}
\section{Introduction}
 The infinite series
\begin{align}
S(r):=\sum_{m=1}^{\infty}\frac{2m}{(m^2+r^2)^2},\quad (r>0),\label{Mathieu-Series}
\end{align}
is called a Mathieu series. It was introduced and studied by \'Emile Leonard Mathieu in his book~\cite{Mathieu} devoted to the elasticity of solid bodies. Since its introduction the series $S(r)$ and its various generalizations have attracted many researchers, who established some remarkable properties of these series including  the various integral representations, the asymptotic expansions, lower and upper estimates, see e.g. Cerone and Lenard~\cite{CL},  Frontczak~\cite{Fro}, Milovanovi\'c and Pog\'any~\cite{MP}, Pog\'any, Srivastava and  Tomovski~\cite{PST}, and references therein.

An integral representation for  the Mathieu series~\eqref{Mathieu-Series} is given by Emersleben~\cite{Eme} as
\begin{align}
S(r)=\frac 1r\int_{0}^{\infty}\frac{x}{e^x-1} \sin(rx)dx.
\end{align}
The integral representation was used by Elbert~\cite{Elb} to derive
the asymptotic expansion of $S(r)$:
\begin{align}
S(r)=\sum_{m=0}^{\infty}(-1)^m\frac{B_{2m}}{r^{2m+2}}
=\frac{1}{r^2}-\frac{1}{6r^4}\pm\cdots,\quad (r\rightarrow\infty),
\end{align}
where $B_{2n}$ denote the even indexed Bernoulli numbers defined by the generating function
\begin{align}
\frac{x}{e^x-1}=\sum_{n=0}^{\infty}B_n\frac{x^n}{n!},\quad |x|<2\pi.
\end{align}

Throughout the paper, we always use notation $\psi(z)=\frac{\Gamma'(z)}{\Gamma(z)}$.
Let $(a_n)_{n\ge 1}$ and $(b_n)_{n\ge 0}$ be two sequences of real ~(or complex) numbers with $a_n\neq 0 $ for all $n\in\mathbb{N}$. The generalized continued fraction
\begin{align}
\tau=b_0+\frac{a_1}{b_1+\frac{a_2}{b_2+\ddots}}=b_0+
\begin{array}{ccccc}
a_1 && a_2 &       \\
\cline{1-1}\cline{3-3}\cline{5-5}
 b_1 & + & b_2 & + \cdots
\end{array}
=b_0+\K_{n=1}^{\infty}
\left(\frac{a_n}{b_n}\right)
\end{align}
is defined as the limit of the $n$th approximant
\begin{align}
\frac{A_n}{B_n}=b_0+\K_{k=1}^{n}\left(\frac{a_k}{b_k}\right)
\end{align}
as $n$ tends to infinity. For the theory of continued fraction, see  Cuyt, Petersen, Verdonk, Waadeland and Jones~\cite{CPV} or Lorentzen and Waadeland~\cite{LW}.

Let $r>0$ and $\Re x>\frac 12$. Let the continued fraction
$CF(r;x)$ with a parameter $r$ be defined by
\begin{align}
CF(r;x)=\frac{1}{\left(x-\frac 12\right)^2+\frac 14\left(1+4r^2\right)+\K_{n=1}^{\infty}
\left(\frac{\kappa_n}{\left(x-\frac 12\right)^2+\lambda_n}\right)},
\label{CF-definition}
\end{align}
where for $n\ge 1$
\begin{align}
\kappa_n= -\frac{n^4\left(n^2 + 4 r^2\right)}{4 (2 n - 1) (2 n + 1)},\quad\lambda_n=\frac 14 (2 n^2 + 2 n + 1 + 4 r^2).
\end{align}
The main purpose of this paper is to establish the following continued fraction expression of the Mathieu series.
\begin{thm} Let $r>0$ and $CF(r;x)$ be defined as~\eqref{CF-definition}. For all positive integer $k$, we have
\begin{align}
S(r)=\sum_{m=1}^{k-1}\frac{2m}{\left(m^2+r^2\right)^2}+CF(r;k),
\label{main-identity}
\end{align}
where the sum for $k-1 = 0$ is stipulated to be zero.
In particular,
\begin{align}
S(r)=CF(r;1).
\end{align}
\end{thm}

\bigskip

\section{Some preliminary lemmas }
In order to prove Theorem 1, we will prepare some lemmas. The following continued fraction formula of Ramanujan plays an important role in the proof of Theorem 1.
\begin{lem} Let $x,l,m$, and $n$ denote complex numbers. Define
\begin{align}
&P=P(x,l,m,n)\\
&=\frac{\Gamma\left(\frac 12(x+l+m+n+1)\right)\Gamma\left(\frac 12(x+l-m-n+1)\right)\Gamma\left(\frac 12(x-l+m-n+1)\right)\Gamma\left(\frac 12(x-l-m+n+1)\right)}{\Gamma\left(\frac 12(x-l-m-n+1)\right)\Gamma\left(\frac 12(x-l+m+n+1)\right)\Gamma\left(\frac 12(x+l-m+n+1)\right)\Gamma\left(\frac 12(x+l+m-n+1)\right)}.\nonumber
\end{align}
Then if either $l,m$, or $n$ is an integer or if $\Re x>0$,
\end{lem}
\begin{align}
\frac{1-P}{1+P}=\frac{2lmn}{x^2-l^2-m^2-n^2+1+\K_{j=1}^{\infty}
\left(
\frac{4(l^2-j^2)(m^2-j^2)(n^2-j^2)}{(2j+1)\left(x^2-l^2-m^2-n^2
+2j^2+2j+1\right)}\right)
}.
\end{align}
\proof This is Entry 35 of B. C. Berndt~\cite{Ber}, p.~157, which was claimed first by Ramanujan~\cite{Ram,Ram1}. The first published proof was provided by Watson~\cite{Wat}. For the full proof of
Entry 35, we refer the reader to L. Lorentzen's paper~\cite{Jac}.\qed

\begin{lem} The canonical contraction of $b_0+\K\left(a_n/b_n\right)$ with
\begin{align*}
C_k=A_{2k},\quad D_k=B_{2k}\quad \mbox{for $k=0,1,2,\ldots$}
\end{align*}
exists if and only if $b_{2k}\neq 0$ for $k=0,1,2,\ldots$, and is given by
\begin{align}
b_0+\K_{n=1}^{\infty}\left(\frac{a_n}{b_n}\right)=&b_0+
\begin{array}{ccccccc}
a_1 b_2& & a_2a_3b_4& & a_4a_5b_2b_6 &\\
\cline{1-1}\cline{3-3}\cline{5-5}\cline{7-7}
a_2+b_1b_2&-& a_3b_4+b_2(a_4+b_3b_4)&-&a_5b_6+b_4(a_6+b_5b_6)&
-\cdots
\end{array}\nonumber\\
&
\begin{array}{ccc}
& a_{2n}a_{2n+1}b_{2n-2}b_{2n+2}& \\
\cline{2-2}
-&a_{2n+1}b_{2n+2}+b_{2n}(a_{2n+2}+b_{2n+1}b_{2n+2})&-\cdots
\end{array}.\nonumber
\end{align}
\end{lem}
\proof It follows from Theorem 12 and Eq.~(2.4.3) of L. Lorentzen, H. Waadeland~\cite{LW}, page 83\nobreakdash--84. For some applications, interested readers may refer to Berndt~[5, p.~121, Eq.~(14.2)] or [5, p.~157]. For a canonical contraction of a continued fraction and the related definitions, also see L. Lorentzen, H. Waadeland~\cite{LW}, page 83.\qed

\begin{lem}
$b_0+\K\left(a_n/b_n\right)\approx d_0+\K\left(c_n/d_n\right)$ if and only if there exists a sequence $\{r_n\}$ of complex numbers with $r_0=1,r_n\neq 0$ for all $n\in \mathbb{N}$, such that
\begin{align}
d_0=b_0,\quad c_n=r_{n-1}r_n a_n,\quad d_n=r_n b_n\quad \mbox{for all
 $n\in \mathbb{N}$.}
\end{align}
\end{lem}
\proof See Theorem 9 of L. Lorentzen, H. Waadeland~\cite{LW}, p.~73.\qed

\section{The proof of Theorem 1}
\begin{lem} Let $r>0$ and $\Re x>\frac 12$, then
\begin{align}
CF(r;x)-CF(r;x+1)
=\frac{2x}{\left(x^2+r^2\right)^2}.
\end{align}
\end{lem}
\begin{rem}
In fact, Lemma 4 was guessed first by the \emph{multiple-correction method} developed in ~\cite{Cao1,CY}. Mortici~\cite{Mor-M} made an important contribution in this direction.
\end{rem}
\proof By applying Lemma 1 with $(x,l)=(2x-1,2r i)$ and dividing both sides by $4rmn i$, we obtain that for $\Re x>\frac 12$
\begin{align}
\frac{1}{4r i}\frac{1-P}{mn(1+P)}
=\frac{1}{(2x-1)^2+4r^2-m^2-n^2+1+\K_{j=1}^{\infty}\left(
\frac{4(-4r^2-j^2)(m^2-j^2)(n^2-j^2)}{(2j+1)
\left((2x-1)^2+4r^2-m^2-n^2
+2j^2+2j+1\right)}\right)}.
\end{align}
Now let $m$ tend to zero and $n$ tend to zero, successively. On the right side, we arrive at
\begin{align}
\frac{1}{(2x-1)^2+1+4r^2+\K_{j=1}^{\infty}\left(
\frac{-4j^4(j^2+4r^2)}{(2j+1)\left((2x-1)^2
+2j^2+2j+1+4r^2\right)}\right)}.\label{R-Hand}
\end{align}
On the other hand, from the definition of $P$, we see easily that $\lim_{m\rightarrow 0}P=1$. A direct calculation with the use of L'Hospital's rule gives
\begin{align}
&\lim_{m\rightarrow 0}\frac{1-P}{m(1+P)}=\lim_{m\rightarrow 0}\frac{1}{1+P}\lim_{m\rightarrow 0}\frac{1-P}{m}=\frac{1}{2} \lim_{m\rightarrow 0}\frac {1-P}{m}
= \frac{1}{2} \lim_{m\rightarrow 0}\frac {\partial}{\partial m}(1-P)  \label{First-Limit}\\
=&\frac{1}{2}\left\{-\psi( -\frac{n}{2}+x-r i)+
\psi(\frac{n}{2}+x-r i)+
\psi(-\frac{n}{2}+x+r i)-
\psi(\frac{n}{2}+x+r i)\right\}.\nonumber
\end{align}
By making use of L'Hospital's rule again, and noting the following classical representation~(e.g., see [1, Eq.~6.3.16, p.~259])
\begin{align}
\psi(z+1)=-\gamma+\sum_{k=1}^{\infty}\left(\frac 1k-\frac{1}{k+z}\right),\quad (z\neq -1,-2,-3,\ldots),
\end{align}
where $\gamma$ denotes Euler-Mascheroni constant, it follows from~\eqref{First-Limit} that
\begin{align}
\lim_{n\rightarrow 0}\lim_{m\rightarrow 0}
\frac{1-P}{mn(1+P)}=&
\lim_{n\rightarrow 0}
\frac {\partial}{\partial n}
\left(\lim_{m\rightarrow 0}\frac 1m\frac{1-P}{1+P}\right)
\label{L-Hand}\\
=&\frac 12\left(\psi'(x-r i) - \psi'(x+r i)\right)\nonumber\\
=&\frac 12 \left(\sum_{k=0}^{\infty}\frac{1}{(x-r i+k)^2}-\sum_{k=0}^{\infty}\frac{1}{(x+r i+k)^2}\right)\nonumber\\
=&2r i\sum_{k=0}^{\infty}\frac{x+k}{\left((x+k)^2+r^2\right)^2}.\nonumber
\end{align}
Combining~\eqref{R-Hand} and ~\eqref{L-Hand}, we get that for $\Re x>\frac 12$
\begin{align}
\sum_{k=0}^{\infty}\frac{2(x+k)}{\left((x+k)^2+r^2\right)^2}=&
\frac{4}{(2x-1)^2+1+4r^2+\K_{j=1}^{\infty}\left(
\frac{-4j^4(j^2+4r^2)}{(2j+1)\left((2x-1)^2
+2j^2+2j+1+4r^2\right)}\right)}\\
=&CF(r;x).\nonumber
\end{align}
Here we used Lemma 3 in the last equality. It is not difficult to check that for $\Re x>\frac 12$
\begin{align}
CF(r;x)-CF(r;x+1)=\sum_{k=0}^{\infty}
\frac{2(x+k)}{\left((x+k)^2+r^2\right)^2}-
\sum_{k=0}^{\infty}\frac{2(x+k+1)}{\left((x+k+1)^2+r^2\right)^2}
=\frac{2x}{\left(x^2+r^2\right)^2}.
\end{align}
This completes the proof of Lemma 4.

\bigskip
\emph{Proof of Theorem 1.} We use the telescoping method. Theorem 1 follows readily from Lemma 4.

\section{Some new inequalities for the Mathieu series}

The bounds for the Mathieu series attracted many mathematicians like
Schr\"oder~\cite{Sch}, Emersleben~\cite{Eme}, Makai~\cite{Mak} and
Diananda~\cite{Dia}. In the past twenty years, many authors like  Alzer, Bagdasaryan, Brenner, Guo, Lampret, Milovanovi\'c,  Mortici, Pog\'any, Qi, Ruehr, Srivastava, Tomovski, etc. have made important contributions to this research topic, see e.g~\cite{ABR,Bag,Lam,MP,Mor-M,PST,PTL} and references therein. Let us briefly recall some simple results.

Mathieu~\cite{Mathieu} himself conjectured only the upper bound $S(r)<r^{-2}$, $r>0$, proved first by Berg~\cite{Berg}. Makai~\cite{Mak} showed the double sided inequalities
\begin{align}
\frac{1}{r^2+\frac 12}<S(r)<\frac{1}{r^2+\frac 16}.
\end{align}
Alzer \emph{et al}.~\cite{ABR} improved the lower bound to
\begin{align}
\frac{1}{r^2+\frac {1}{2\zeta(3)}}<S(r)<\frac{1}{r^2+\frac 16},
\label{Alzer's result}
\end{align}
where the constant $1/(2\zeta(3))$ and $1/6$ are sharp.

Milovanovi\'c and Pog\'any~\cite{MP} stated a composite upper bound
of simple structure,
\begin{align}
S(r)\le
\begin{cases}
\frac{1}{r^2+\frac 14},&0\le r\le \frac{\sqrt 3}{2},\\
\frac{1}{\sqrt{1+4r^2}-1},&r>\frac{\sqrt 3}{2},
\end{cases}
\end{align}
which is superior to \eqref{Alzer's result} in the interval $[0,\sqrt{(5+2\sqrt{3})/6}\approx 1.18772)$.

Let $a_1=1$, $b_1=\left(x-\frac 12\right)^2-\frac 14 +r^2$, for $n\ge 1$
\begin{align}
a_{2n+1}=\frac{n\left(n^2 + 4 r^2\right)}{2(2 n + 1)},\quad a_{2n}=\frac{n^3}{2 (2 n - 1) },\label{a(n)-def}
\end{align}
and
\begin{align}
b_{2n+1}=\left(x-\frac 12\right)^2-\frac 14 +\frac{r^2}{2 n+1},
\quad b_{2n}=1.\label{b(n)-def}
\end{align}

By Lemma 2, it is not difficult to prove that
\begin{align}
CF(r;x)=\K_{n=1}^{\infty}\left(\frac{a_n}{b_n}\right).
\end{align}

We let $z=(x-\frac 12)^2-\frac 14$, $c_1=2$, $d_1=2z+2r^2$,
\begin{align}\begin{cases}
&c_{2n+1}=n(n^2+4r^2),\quad c_{2n}=n^3,\\
&d_{2n+1}=2(2n+1)z+2r^2,\quad d_{2n}=1.
\end{cases}
\end{align}

It follows easily from Lemma 3 that

\begin{lem} Let $\Re x>\frac 12$. With the above notation, we have
\begin{align}
CF(r;x)=\K_{n=1}^{\infty}\left(\frac{a_n}{b_n}\right)
=\K_{n=1}^{\infty}\left(\frac{c_n}{d_n}\right).
\end{align}
\end{lem}
\begin{lem} Assume $x>\frac 12$. For all positive integer $l$, then \begin{align}
\K_{n=1}^{2l}\left(\frac{a_n}{b_n}\right)<
CF(r;x)
<\K_{n=1}^{2l-1}\left(\frac{a_n}{b_n}\right).
\end{align}
\end{lem}
\proof As the partial coefficients of the continued fraction $\K_{n=1}^{\infty}\left(\frac{a_n}{b_n}\right)$ are positive, the assertion is deduced readily from the theory of the continued fraction and the first equality in Lemma 5.\qed

The following theorem tells us how to obtain sharp bounds for the Mathieu series.
\begin{thm} Let $r>0$ and $k, l\in \mathbb{N}$. Let two sequences $\left(a(n)\right)_{n\ge 1}, \left(b(n)\right)_{n\ge 1}$ be defined by~\eqref{a(n)-def} and~\eqref{b(n)-def} with $x=k$, respectively. Then
\begin{align}
\sum_{m=1}^{k-1}\frac{2m}{\left(m^2+r^2\right)^2}+
\K_{n=1}^{2l}\left(\frac{a_n}{b_n}\right)<
S(r)<\sum_{m=1}^{k-1}\frac{2m}{\left(m^2+r^2\right)^2}+
\K_{n=1}^{2l-1}\left(\frac{a_n}{b_n}\right).\label{Mathieu-bounds}
\end{align}
In particular,
\begin{align}
&\frac{2}{(1 + r^2)^2}+\frac{1}{5/2 + r^2}<S(r)<\frac{2}{(1 + r^2)^2} +\frac{1}{2 + r^2},\label{Mathieu-bounds-1}\\
&\frac{2}{(1 + r^2)^2}+\frac{4}{(4 + r^2)^2}+\frac{1}{13/2 + r^2}<S(r)<\frac{2}{(1 + r^2)^2}+\frac{4}{(4 + r^2)^2} +\frac{1}{6 + r^2}.\label{Mathieu-bounds-2}
\end{align}
\end{thm}
\proof \eqref{Mathieu-bounds} follows readily  from Theorem 1 and Lemma 6. Taking $(k,l)=(2,1)$ and $(k,l)=(3,1)$ in \eqref{Mathieu-bounds}, respectively, we can obtain \eqref{Mathieu-bounds-1} and \eqref{Mathieu-bounds-2}.\qed

\begin{rem} For comparison, our upper bound in \eqref{Mathieu-bounds-1} improves (4.3) when
$0\le r<\sqrt{-2 + \sqrt{7}}\approx 0.803587$. It is not hard to check that if
\begin{align*}
r\in \left(\sqrt{\frac{-6+5 \zeta(3)}{
2 + \sqrt{-2 + 11 \zeta(3) - 5 \zeta^2(3)}}},\sqrt{
\frac{2+\sqrt{-2 + 11 \zeta(3)- 5 \zeta^2(3)}}{\zeta(3)-1}
}\right)\approx (0.0507096,4.44903),
\end{align*}
then our lower bound in \eqref{Mathieu-bounds-1} is  superior to Alzer's in \eqref{Alzer's result}. In addition,  the bounds in \eqref{Mathieu-bounds-2} are always superior to the bounds in \eqref{Mathieu-bounds-1} for all $r>0$.
\end{rem}

\begin{rem} Taking $k=2$ in Theorem 1, letting $r$ tend to zero, and then applying the second equality in Lemma 5, we can deduce the following continued fraction for Ap\'ery number $\zeta(3)$
\begin{align*}
\zeta(3)=1+
\begin{array}{ccccccccccc}
1&&1^3 & & 1^3 & & 2^3 & & 2^3 &\\
\cline{1-1}\cline{3-3}\cline{5-5}\cline{7-7}\cline{9-9}\cline{11-11}
2^2\cdot 1&+&1 & + & 2^2\cdot 3 & + & 1 & +&2^2\cdot 5 & +\cdots
\end{array}.
\end{align*}
Also see Berndt~\cite{Ber}, p.~155.
\end{rem}

\end{document}